\documentclass[11pt,leqno,fleqn]{article}
\usepackage{amssymb}
\newcommand{\rf}[1]{{\rm (\ref{#1})}}
\newcommand{\dps}{\displaystyle}
\newcommand{\kies}[2]{\mbox{${{#1}\choose{#2}}$}}
\newcommand{\PP}{{\cal P}}
\newcounter{stelling}
\newcommand{\thmz}[1]{\refstepcounter{stelling}\vspace{4mm}\noindent{\bf Theorem \thestelling.}{\it #1}}
\newcommand{\thmnmz}[2]{\refstepcounter{stelling}\vspace{4mm}\noindent{\bf Theorem \thestelling\ {\rm (#1)}.}{\it #2}}
\newcommand{\de}[2]{\begin{equation}\label{#1} #2 \end{equation}}
\newcommand{\dy}[2]{%
\refstepcounter{equation}%
\label{#1}%
\begin{list}{}{
\topsep 3mm
\leftmargin 18mm
\rightmargin 0cm
\itemsep 0mm
\listparindent 0mm
\parsep 0mm
\itemsep 0mm
\labelsep 0mm
\labelwidth 18mm
}%
\item[\rm (\theequation)\hfill]
#2
\end{list}%
}
\newcommand{\dyz}[1]{%
\refstepcounter{equation}%
\begin{list}{}{
\topsep 3mm
\leftmargin 18mm
\rightmargin 0cm
\itemsep 0mm
\listparindent 0mm
\parsep 0mm
\itemsep 0mm
\labelsep 0mm
\labelwidth 18mm
}%
\item[\rm (\theequation)\hfill]
#1
\end{list}%
}
\newcommand{\dyyz}[1]{\dyz{\raggedright$\dps#1$}}
\newcommand{\dyy}[2]{\dy{#1}{\raggedright$\dps#2$}}
\newcounter{sectie}
\newcommand{\sectz}[1]{\refstepcounter{sectie}
\section*{\boldmath \thesectie. #1}%
}
\newcommand{\pf}{\vspace{3mm}\noindent{\bf Proof.}\ }
\newcommand{\bx}{\hspace*{\fill} \hbox{\hskip 1pt \vrule width 4pt height 8pt depth 1.5pt \hskip 1pt}

\addvspace{4mm}}

\newcommand{\oN}{{\mathbb{N}}}
\begin{document}

\begin{center}
{\LARGE On the Shannon capacity of sums and products of graphs}

\bigskip
Alexander Schrijver\footnote{Centrum Wiskunde \& Informatica and University of Amsterdam}

\end{center}

\noindent
{\small
{\bf Abstract}
Let $\Theta(G)$ denote the Shannon capacity of a graph $G$.
We give an elementary proof of the equivalence, for any graphs $G$ and $H$, of
the inequalities $\Theta(G\sqcup H)>\Theta(G)+\Theta(H)$ and
$\Theta(G\boxtimes H)>\Theta(G)\Theta(H)$.
This was shown independently by Wigderson and Zuiddam [2022]  using Kadison-Dubois
duality and the Axiom of choice.
}

\sectz{Introduction}

Let $G$ be a graph.
(All graphs in this paper are undirected and loopless.)
A {\em stable set} in $G$ is a set of pairwise nonadjacent vertices.
The {\em stable set number} $\alpha(G)$ is the maximum cardinality of a stable set in $G$.

The {\em sum} $G+H$ of graphs $G$ and $H$ is the disjoint union of $G$ and $H$.
Trivially,
\de{11ap22g}{
\alpha(G+H)=\alpha(G)+\alpha(H).
}
The {\em strong product} $GH$ of $G$ and $H$ is the graph with vertex set $V(G)\times V(H)$
where distinct $(u,v)$ and $(u',v')$ in $V(G)\times V(H)$ are adjacent if and only if
$u$ and $u'$ are equal or adjacent in $G$ and $v$ and $v'$ are equal or adjacent in $H$.

Since sum and strong product are associative, commutative, and distributive,
this makes the set of graphs to a commutative semiring,
with unit the one-vertex graph $K_1$.
Sum and strong product are often denoted by $G\sqcup H$ and $G\boxtimes H$,
but the semiring notation $G+H$ and $GH$ is more efficient here.

As the cartesian product of stable sets in $G$ and $H$ is a stable set in $GH$ we have
\de{11ap22a}{
\alpha(GH)\geq\alpha(G)\alpha(H),
}
but strict inequality may occur, even if $G=H$ (for instance for $G=H=C_5$, the five-cycle).
This made Shannon [1956] define what is now called the {\em Shannon capacity} $\Theta(G)$
of a graph $G$:
\de{11ap22c}{
\Theta(G):=\sup_{k\in\oN}\alpha(G^k)^{1/k}=\lim_{k\to\infty}\alpha(G^k)^{1/k}.
}
The second equality in \rf{11ap22c} follows from \rf{11ap22a} and Fekete's lemma [1923].
(In fact, Shannon introduced $\log\Theta(G)$ as the `zero-error capacity' of
the `channel' $G$.)

The inequality \rf{11ap22a} implies
\de{11ap22b}{
\Theta(GH)\geq\Theta(G)\Theta(H).
}
Haemers [1979] (disproving a conjecture of Shannon [1956])
gave examples of graphs $G,H$ with strict inequality in \rf{11ap22b}.
(In fact, Haemers showed that the `Schl\"afli graph' $G$ satisfies
$\Theta(G)\Theta(\overline G)<|V(G)|\leq\alpha(G\overline G)$.
Here $\overline G$ is the graph complementary to $G$.)

On the other hand, for each graph $G$ and $n\in\oN$:
\de{12ap22a}{
\Theta(G^n)=\Theta(G)^n,
}
as follows directly from definition \rf{11ap22c}.

The value of $\Theta(C_5)$ was for a long an open question, until Lov\'asz [1979]
introduced the upper bound $\vartheta(G)$ on $\Theta(G)$ yielding $\Theta(C_5)=\sqrt{5}$.
Since $\vartheta(GH)=\vartheta(G)\vartheta(H)$ for all $G,H$, the Haemers examples imply
that $\Theta(G)<\vartheta(G)$ may occur.

As for the sum, Shannon showed that for all graphs $G$ and $H$ one has
\de{11ap22d}{
\Theta(G+H)\geq\Theta(G)+\Theta(H).
}
(For completeness, we give a proof below.)
Shannon conjectured that for all $G,H$ equality holds in \rf{11ap22d}.
This was disproved by Alon [1998], by displaying graphs $G$ and $H$ with
$\Theta(G+H)>\Theta(G)+\Theta(H)$.
In fact, strict inequality holds for any $G$ and $H$ that satisfy
$\Theta(GH)>\Theta(G)\Theta(H)$, as follows (using \rf{12ap22a} and \rf{11ap22d}) from
\dyy{11ap22h}{
\Theta(G+H)^2
=
\Theta((G+H)^2)
=
\Theta(G^2+2GH+H^2)
\geq
\Theta(G^2)+2\Theta(GH)+\Theta(H^2)
=
\Theta(G)^2+2\Theta(GH)+\Theta(H)^2
>
\Theta(G)^2+2\Theta(G)\Theta(H)+\Theta(H)^2
=
(\Theta(G)+\Theta(H))^2.
}
So Haemers' counterexamples $G,H$ for products also work for sums.

In this paper we give an elementary proof of the fact that for all $G,H$:
\dyyz{
\Theta(GH)>\Theta(G)\Theta(H)$~~~$\iff$~~~$\Theta(G+H)>\Theta(G)+\Theta(H).
}
This was proved (independently) by Wigderson and Zuiddam [2022], using Strassen's theory
of asymptotic spectra (based on Kadison-Dubois duality)
and the Axiom of choice.

More strongly, consider any $n\in\oN$ and graphs $G_1,\ldots, G_n$.
Then for any polynomial $p\in\oN[x_1,\ldots,x_n]$ one has
\de{13ap22a}{
\Theta(p(G_1,\ldots,G_n))\geq p(\Theta(G_1),\ldots,\Theta(G_n)).
}
(This follows from \rf{11ap22d} and \rf{11ap22b}.)
Now if equality holds in \rf{13ap22a} for one polynomial $p$ in which each of the variables
$x_1,\ldots,x_n$ occurs, then equality holds in \rf{13ap22a} for all polynomials $p$.

\sectz{Shannon's inequality}

For self-containedness of this paper, we give a proof of Shannon's inequality:

\thmnmz{Shannon [1956]}{
$\Theta(G+H)\geq \Theta(G)+\Theta(H)$.
}

\pf
For all $n,t\geq 1$, using \rf{11ap22g} and \rf{11ap22a}:
\dyy{11ap22f}{
\alpha((G+H)^n)
=
\alpha(\sum_{k=0}^n\kies{n}{k}G^kH^{n-k})
=
\sum_{k=0}^n\kies{n}{k}\alpha(G^kH^{n-k})
\geq
\sum_{k=0}^n\kies{n}{k}\alpha(G^k)\alpha(H^{n-k})
\geq
\sum_{k=0}^n\kies{n}{k}\alpha(G^t)^{\lfloor k/t\rfloor}\alpha(H^t)^{\lfloor (n-k)/t\rfloor}
\geq
\sum_{k=0}^n\kies{n}{k}
\alpha(G^t)^{k/t}\alpha(G^t)^{-1}
\alpha(H^t)^{(n-k)/t}\alpha(H^t)^{-1}
=
(\alpha(G^t)^{1/t}
+\alpha(H^t)^{1/t})^n
\alpha(G^t)^{-1}\alpha(H^t)^{-1}.
}
So for each $t\geq 1$:
\dyyz{
\Theta(G+H)
=\sup_{n\in\oN}\alpha((G+H)^n)^{1/n}
\geq
\sup_{n\in\oN}
(\alpha(G^t)^{1/t}+\alpha(H^t)^{1/t})\alpha(G^t)^{-1/n}\alpha(H^t)^{-1/n}
=
\alpha(G^t)^{1/t}+\alpha(H^t)^{1/t}.
}
So $t\to\infty$ gives the theorem.
\bx

\sectz{
$\Theta(GH)>\Theta(G)\Theta(H)$~~$\iff$~~$\Theta(G+H)>\Theta(G)+\Theta(H)$
}

\thmz{
$\Theta(GH)>\Theta(G)\Theta(H)$ if and only if $\Theta(G+H)>\Theta(G)+\Theta(H)$.
}

\pf
Necessity follows from \rf{11ap22h}
To see sufficiency, assume $\Theta(GH)\leq\Theta(G)\Theta(H)$.
Then for all $i,j\in\oN$, using \rf{11ap22b} and \rf{12ap22a}:
\dyyz{
\Theta(G^iH^j)\Theta(G)^j\Theta(H)^i
=
\Theta(G^iH^j)\Theta(G^j)\Theta(H^i)\leq\Theta((GH)^{i+j})
=\Theta(GH)^{i+j}\leq\Theta(G)^{i+j}\Theta(H)^{i+j}.
}
So $\Theta(G^iH^{j})\leq\Theta(G)^i\Theta(H)^{j}$.
Hence for each $n$, using \rf{11ap22g}:
\dyyz{
\alpha((G+H)^n)
=
\alpha(\sum_{k=0}^n\kies{n}{k}G^kH^{n-k})
=
\sum_{k=0}^n\kies{n}{k}\alpha(G^kH^{n-k})
\leq
\sum_{k=0}^n\kies{n}{k}\Theta(G^kH^{n-k})
\leq
\sum_{k=0}^n\kies{n}{k}\Theta(G)^k\Theta(H)^{n-k}
=
(\Theta(G)+\Theta(H))^n.
}
Taking $n$-th roots and $n\to\infty$ gives $\Theta(G+H)\leq\Theta(G)+\Theta(H)$.
\bx

\sectz{Extension to polynomials}

We also give an elementary proof of the following extension that was shown by
Wigderson and Zuiddam [2022] using Kadison-Dubois duality and the Axiom of choice.

For given graphs $G_1,\ldots,G_n$, define
\de{4fe18b}{
\PP=\{p\in\oN[x_1,\ldots,x_n]\mid \Theta(p(G_1,\ldots,G_n))=p(\Theta(G_1),\ldots,\Theta(G_n))\}.
}

\thmz{
Let $G_1,\ldots,G_n$ be graphs with at least one vertex.
Then
$\PP=\oN[x_1,\ldots,x_n]$
if and only if
$\PP$ contains a polynomial in which all variables $x_1,\ldots,x_n$ occur.
}

\pf
Necessity being trivial, we prove sufficiency.
Let $\underline G:=(G_1,\ldots,G_n)$ and $\Theta(\underline G):=(\Theta(G_1),\ldots,\Theta(G_n))$.
So $p(\Theta(\underline G))\leq \Theta(p(\underline G))$ for any polynomial $p\in\oN[x_1,\ldots,x_n]$.

We first show that for $p,q\in\oN[x_1,\ldots,x_n]$:
\dy{4fe18bx}{
if $p+q\in\PP$, then $p\in\PP$.
}
Indeed,
\dyyz{
\Theta((p+q)(\underline G))
=
(p+q)(\Theta(\underline G))
=
p(\Theta(\underline G))+q(\Theta(\underline G))
\leq
\Theta(p(\underline G))+\Theta(q(\underline G))
\leq
\Theta(p(\underline G)+q(\underline G))
=
\Theta((p+q)(\underline G)).
}
Hence we have equality throughout, implying $\Theta(p(\underline G))=p(\Theta(\underline G))$.
This proves \rf{4fe18bx}.

Similarly,
\dy{4fe18by}{
if $pq\in\PP$ and $q\neq 0$, then $p\in\PP$.
}
Indeed,
\dyyz{
\Theta((pq)(\underline G))
=
(pq)(\Theta(\underline G))
=
p(\Theta(\underline G))q(\Theta(\underline G))
\leq
\Theta(p(\underline G))\Theta(q(\underline G))
\leq
\Theta(p(\underline G)q(\underline G))
=
\Theta((pq)(\underline G)).
}
Hence we have equality throughout, implying $\Theta(p(\underline G))=p(\Theta(\underline G))$.
This proves \rf{4fe18by}.

Moreover, for $p\in\oN[x_1,\ldots,x_n]$ and $k\in\oN$,
\dy{4fe18c}{
if $p\in\PP$ then $p^k\in\PP$.
}
Indeed, if $p\in\PP$, then
\dyyz{
\Theta(p^k(\underline G))
=
\Theta(p(\underline G)^k)
=
(\Theta(p(\underline G)))^k
=
(p(\Theta(\underline G)))^k
=
(p^k(\Theta(\underline G))),
}
proving \rf{4fe18c}.

Now let $p\in\PP$ with each $x_1,\ldots,x_n$ occurring in $p$.
Then for some $k\in\oN$, $p^k$ contains as term a monomial $q$ in which
each variable occurs at least once.
As $p^k\in\PP$ by \rf{4fe18c}, we know by \rf{4fe18bx} that $q\in\PP$.
Now for each monomial $\mu$ in $\oN[x_1,\ldots,x_n]$ there exists a large
enough $N$ such that $\mu$ is a divisor of $q^N$.
So by \rf{4fe18b}, each monomial belongs to $\PP$.

Now consider any polynomial $r=q_1+\cdots+q_t$ in $\oN[x_1,\ldots,x_n]$,
where each $q_i$ is a monomial.
Then for each $i_1,\ldots,i_t\in\oN$, $\mu:=\prod_{j=1}^tq_j^{i_j}$ is a monomial, implying
\dyyz{
\Theta(\prod_{j=1}^tq_j(\underline G)^{i_j})
=
\Theta(\mu(\underline G))
=
\mu(\Theta(\underline G))
=
\prod_{j=1}^tq_j(\Theta(\underline G))^{i_j}.
}
This implies, for each $k\in\oN$, using the additivity (\rf{11ap22g})
of the function $\alpha$:
\dyyz{
\alpha(r(\underline G)^k)
=
\alpha((\sum_{j=1}^tq_j(\underline G))^k)
=
\alpha(
\sum_{i_1,\ldots,i_t\in\oN\atop i_1+\cdots+i_t=k}
\kies{k}{i_1,\ldots,i_t}
\prod_{j=1}^tq_j(\underline G)^{i_j})
=
\sum_{i_1,\ldots,i_t\in\oN\atop i_1+\cdots+i_t=k}
\kies{k}{i_1,\ldots,i_t}
\alpha(\prod_{j=1}^tq_j(\underline G)^{i_j})
\leq
\sum_{i_1,\ldots,i_n\in\oN\atop i_1+\cdots+i_t=k}
\kies{k}{i_1,\ldots,i_t}
\Theta(\prod_{j=1}^tq_j(\underline G)^{i_j})
=
\sum_{i_1,\ldots,i_n\in\oN\atop i_1+\cdots+i_t=k}
\kies{k}{i_1,\ldots,i_t}
\prod_{j=1}^tq_j(\Theta(\underline G))^{i_j}
=
\big(\sum_{j=1}^tq_j(\Theta(\underline G))\big)^k
=
(r(\Theta(\underline G)))^k.
}
Taking $k$-th roots and $k\to\infty$ we obtain
$\Theta(r(\underline G))\leq r(\Theta(\underline G))$.
So $r\in\PP$.
\bx

\section*{References}\label{REF}
\markboth{References}{References}
{\small
\begin{itemize}{}{
\setlength{\labelwidth}{4mm}
\setlength{\parsep}{0mm}
\setlength{\itemsep}{0mm}
\setlength{\leftmargin}{5mm}
\setlength{\labelsep}{1mm}
}
\item[\mbox{\rm[1998]}] N. Alon, 
The Shannon capacity of a union,
{\em Combinatorica} 18 (1998) 301--310.

\item[\mbox{\rm[1923]}] M. Fekete,
\"Uber die Verteilung der Wurzeln bei gewissen algebraischen
Gleichungen mit ganzzahligen Koeffizienten,
{\em Mathematische Zeitschrift} 17 (1923) 228--249.

\item[\mbox{\rm[1979]}] W. Haemers, 
On some problems of Lov\'asz concerning the Shannon capacity of a graph,
{\em {IEEE} Transactions on Information Theory} {IT}-25 (1979) 231--232.

\item[\mbox{\rm[1979]}] L. Lov\'asz, 
On the Shannon capacity of a graph,
{\em {IEEE} Transactions on Information Theory} {IT}-25 (1979) 1--7.

\item[\mbox{\rm[1956]}] C.E. Shannon, 
The  zero error capacity of a noisy channel,
{\em {IRE} Transactions on Information Theory} {IT}-2 (1956) 8--19.

\item[\mbox{\rm[2022]}] A. Wigderson and J. Zuiddam, 
Asymptotic spectra: theory, applications and extensions,
manuscript, 2022.

\end{itemize}
}

\end{document}